\documentclass{article}
\usepackage[utf8]{inputenc}
\usepackage[utf8]{inputenc}
\usepackage[utf8]{inputenc}
\usepackage{dirtytalk}
\usepackage{tikz}
\usepackage{geometry}
\usepackage{graphicx}
\usepackage{enumitem}
\usepackage{parskip}
\usepackage{amsfonts}
\usepackage{amssymb}
\usepackage{amsmath}
\usepackage{amsthm}
\usepackage{xcolor}
\usepackage{hyperref}
\newlist{myitemize}{itemize}{1}
\setlist[myitemize,1]{leftmargin = 0.5in}
\hypersetup{colorlinks = true, linkcolor = red,  linktocpage}

\theoremstyle{plain}
\newtheorem{thm}{Theorem}[section]
\newtheorem*{thm*}{Theorem}

\newtheorem{cor}[thm]{Corollary}

\theoremstyle{definition}

\newtheorem{conj}[thm]{Conjecture}
\newtheorem{exmp}[thm]{Example}
\newtheorem{rem}[thm]{Remark}

\title{\textbf{\small{COUNTING THE NUMBER OF INTEGRAL FIXED POINTS OF A DISCRETE DYNAMICAL SYSTEM WITH APPLICATIONS FROM ARITHMETIC STATISTICS, I}}}
\author{\footnotesize{BRIAN KINTU}}
\date{\small{\textit{Happily Dedicated: Dept. of MCS \& Presidency of Meric S. Gertler at the University of Toronto}}}

\begin{document}
\maketitle
\begin{abstract}
\small{In this article, we inspect a surprising relationship between the set of fixed points of a polynomial map $\varphi_{d, c}$ defined by $\varphi_{d, c}(z) = z^d + c$ for all $c, z \in \mathbb{Z}$ and the coefficient $c$, where $d > 2$ is an integer. Inspired greatly by the elegant counting problems along with the very striking results of Bhargava-Shankar-Tsimerman and their collaborators in arithmetic statistics, and also by interesting point-counting result of Narkiewicz on rational periodic points of any odd degree map $\varphi_{d, c}$ in arithmetic dynamics, we then first prove that for any prime $p\geq 3$, the average number of distinct integral fixed points of any $\varphi_{p, c}$ modulo $p$ is $3$ or $0$ as $c$ tends to infinity. Inspired further by a conjecture of Hutz on rational periodic points of $\varphi_{p-1, c}$ for any prime $p\geq 5$ in arithmetic dynamics, we then also prove that the average number of distinct integral fixed points of any $\varphi_{p-1, c}$ modulo $p$ is $1$ or $2$ or $0$ as $c\to \infty$. Finally, we then apply here density and number field-counting results from arithmetic statistics, and as a result obtain counting and statistical results on the irreducible monic integer polynomials and number fields arising naturally in our polynomial discrete dynamical settings.}
\end{abstract}

\begin{center}
\tableofcontents
\end{center}
\begin{center}
    \section{Introduction}
\end{center}
\noindent
Consider any morphism $\varphi: {\mathbb{P}^N(K)} \rightarrow {\mathbb{P}^N(K)} $ of degree $d \geq 2$ defined on a projective space ${\mathbb{P}^N(K)}$ of dimension $N$, where $K$ is a number field. Then for any $n\in\mathbb{Z}$ and $\alpha\in\mathbb{P}^N(K)$, we call $\varphi^n = \underbrace{\varphi \circ \varphi \circ \cdots \circ \varphi}_\text{$n$ times}$ the $n^{th}$ \textit{iterate of $\varphi$} and call $\varphi^n(\alpha)$ the \textit{$n^{th}$ iteration of $\varphi$ on $\alpha$}. By convention, $\varphi^{0}$ acts as the identity map, i.e., $\varphi^{0}(\alpha) = \alpha$ for every point $\alpha\in {\mathbb{P}^N(K)}$. The everyday philosopher may want to know (quoting here Devaney \cite{Dev}): \say{\textit{Where do points $\alpha, \varphi(\alpha), \varphi^2(\alpha), \ \cdots\ ,\varphi^n(\alpha)$ go as $n$ becomes large, and what do they do when they get there?}} So now, for any given integer $n\geq 0$ and for any given point $\alpha\in {\mathbb{P}^N(K)}$, we then call the set consisting of all the iterates $\varphi^n(\alpha)$ the \textit{(forward) orbit of $\alpha$}; and which in the theory of dynamical systems we usually denote this set by:
\begin{equation}
\mathcal{O}^{+}(\alpha):= \{\varphi^n(\alpha) : n \in \mathbb{Z}_{\geq 0} \}.
\end{equation}
\indent One of the principal goals in the area of \say{arithmetic dynamics} (a newly emerging area of mathematics concerned with studying number-theoretic properties of discrete dynamical systems) is to classify all the points $\alpha\in\mathbb{P}^N(K)$ according to the behavior of their forward orbits $\mathcal{O}^{+}(\alpha)$. In this direction, we call any point $\alpha\in {\mathbb{P}^N(K)}$ a \textit{periodic point of $\varphi$}, whenever $\varphi^n (\alpha) = \alpha$ for some $n\in \mathbb{Z}_{\geq 0}$; and call any integer $n\geq 0$ such that $\varphi^n (\alpha) = \alpha$ a \textit{period of $\alpha$}, and the smallest such positive integer $n\geq 1$ is then called the \textit{exact period of $\alpha$}. We denote the set of all periodic points of $\varphi$ by Per$(\varphi, {\mathbb{P}^N(K)})$; and for any given point $\alpha\in$Per$(\varphi, {\mathbb{P}^N(K)})$ we then call the set of all iterates of $\varphi$ on $\alpha$, \textit{a periodic orbit of $\alpha$}. In Section \ref{sec4}, \ref{sec5}, \ref{sec6}, \ref{sec7} and \ref{sec8}, we study counting questions that are greatly inspired by all the beautiful work of Bhargava-Shankar-Tsimerman (BST) and their collaborators in the area of \say{arithmetic statistics} (an area of mathematics concerned with studying the distributions of arithmetic objects (quantities)); and among such questions includes the natural question: \say{\textit{How many distinct integral fixed orbits can any $\varphi_{p,c}$ and $\varphi_{p-1,c}$ acting independently on the space $\mathbb{Z} / p\mathbb{Z}$ via iteration have on average as $c\to \infty$?}} In doing so, we then first prove that conditioning on Narkiewicz's Theorem \ref{theorem 3.2.1} yields the following main theorem on every map $\varphi_{p,c}$; stated later more precisely as Theorem \ref{2.2}:

\begin{thm}\label{Binder-Brian1}
Let $p\geq 3$ be any fixed prime integer, and assume Theorem \ref{theorem 3.2.1}. Let $\varphi_{p, c}$ be a map defined by $\varphi_{p, c}(z) = z^p+c$ for all $c, z\in\mathbb{Z}$. The number of distinct integral fixed points of any $\varphi_{p,c}$ modulo $p$ is $3$ or zero.
\end{thm}

Inspired further by the activity and down-to-earth work of (BST) in arithmetic statistics, and also by an intriguing conjecture of Hutz \ref{conjecture 3.2.1} (though not more importantly in our case attempting to prove his Conjecture \ref{conjecture 3.2.1}) along with recent work of Panraksa \cite{par2} on rational periodic points of polynomial maps $\varphi_{d, c}$ of any even degree $d> 2$ in arithmetic dynamics, we then revisit the setting in Section \ref{sec2} and then consider in Section \ref{sec3} any polynomial map $\varphi_{p-1,c}$ iterated on the space $\mathbb{Z}\slash p\mathbb{Z}$. In doing so, we then also prove unconditionally the following main theorem on every polynomial map $\varphi_{p-1,c}$; which we state later more precisely as Theorem \ref{6.0.2}:

\begin{thm}\label{Binder-Brian2}
Let $p\geq 5$ be any fixed prime integer, and let $\varphi_{p-1, c}$ be a map defined by $\varphi_{p-1, c}(z) = z^{p-1}+c$ for all $c, z\in\mathbb{Z}$. Then the number of distinct integral fixed points of any map $\varphi_{p-1,c}$ modulo $p$ is $1$ or $2$ or zero.
\end{thm}

\noindent Notice that the count obtained in Theorem \ref{Binder-Brian2} on the number of distinct integral fixed points of any $\varphi_{p-1,c}$ modulo $p$ is independent of $p$ (and so independent of degree of $\varphi_{p-1,c})$ in each of the three possibilities. Moreover, we may also observe that the expected total count (namely, $1+2+0 =3$) in Theorem \ref{Binder-Brian2} on the number of distinct integral fixed points in the whole family of maps $\varphi_{p-1,c}$ modulo $p$ is also independent of $p$ and deg$(\varphi_{p-1,c})$; an observation which somewhat surprisingly coincides not only with a similar observation on the count in each of the two possibilities in Theorem \ref{Binder-Brian1} but also coincides with a similar observation on the expected total count (namely, $3+0 =3$) on the number of distinct integral fixed points in the whole family of maps $\varphi_{p,c}$ modulo $p$. 

Since we know the inclusion $\mathbb{Z}\hookrightarrow\mathbb{Z}_{p}$ of rings, and so the space $\mathbb{Z}_{p}$ of all $p$-adic integers is evidently a much larger space than $\mathbb{Z}$. So then, inspired by the work of Adam-Fares \cite{Ada} in arithmetic dynamics and again by a\say{counting-application} philosophy in arithmetic statistics, we again inspect in a paper \cite{BK3} the aforementioned relationship where it is $\mathbb{Z}_{p}$ that's considered. Interestingly, we again obtain the same counting and asymptotics in the setting when $\varphi_{p-1,c}$ is iterated on $\mathbb{Z}_{p}\slash p\mathbb{Z}_{p}$, and get very different counting and asymptotics in the case when $\varphi_{p,c}$ iterated on $\mathbb{Z}_{p}\slash p\mathbb{Z}_{p}$. Motivated by a $K$-rational periodic point-counting result of Narkiewicz \cite{Narkie} on maps $\varphi_{p,c}$ defined over any real algebraic number field $K$ of degree $n\geq 2$, in a forthcoming work \cite{BK2}, we revisit the counting setting in Section \ref{sec2} and then consider any $\varphi_{p^{\ell},c}$ defined over any real algebraic number $K\slash \mathbb{Q}$ of any degree $n\geq 2$ for any prime $p\geq 3$ and any integer $\ell \geq 1$. In doing so, we show that conditioning on Narkiewicz' theorem and again using the same elementary counting technique, we can obtain a fixed integral point-counting result that is not only independent of the degree $n$ of any real $K$ and any degree $p^{\ell}$ of any map $\varphi_{p^{\ell},c}$, but is also very analogous to Theorem \ref{Binder-Brian1} for any real number field $K$, any odd prime $p$ and any integer $\ell$. More in that work \cite{BK2}, we again revisit the counting setting in Section \ref{sec3} and then consider any polynomial map $\varphi_{(p-1)^{\ell},c}$ defined over any algebraic number $K\slash \mathbb{Q}$ (where $K$ needn't be a real algebraic number field) of any degree $n\geq 2$, for any prime integer $p\geq 5$ and any integer $\ell \geq 1$. In doing so, we find that we can again obtain a fixed integral point-counting result that's not only independent of both the degree $n$ of any number field $K$ and any prime $p$, but is also very similar to Theorem \ref{Binder-Brian2} for any given algebraic number field $K$, any prime $p$ and any integer $\ell$.

It's worth mentioning that there are several authors in the literature who have also done an explicit study on orbits; and among such authors includes, Silverman \cite{Sil} on the number of integral points in forward orbits of rational functions of degree $\geq 2$ defined over the field $\mathbb{C}$, and Wade \cite{Wade} on the average number of integral points in forward orbits of rational functions of degree $\geq 2$ defined over $K$. In his 1993 beautiful paper \cite{Sil}, Silverman showed [\cite{Sil}, Theorem A] that the forward orbit of a rational function of degree $\geq 2$ whose second iterate is not a polynomial over $\mathbb{C}$, contains only finitely many integer points. Motivated by a point-counting result of Silverman \cite{Sil} and conditioning on a standard height uniform conjecture in arithmetic geometry, Wade \cite{Wade} established a zero-average result on the number of integral points in the forward orbit of a rational function of degree $\geq 2$ over $K$. Now observe that in each of these works \cite{Sil, Wade} the focus is on understanding very well integral interiority of forward orbits. In our case, the focus is on understanding a somewhat less thorny problem, namely, as we recall from the above: the problem of counting the number of distinct integral fixed orbits (on average), and then hopefully understand the statistical behavior (along with the question of the meaning) of the achieved count. Such a problem may supposedly be of some interesting insight in the area of classical dynamical systems, since one of the main objectives in that area is to understand \textit{all} orbits via topological and analytic techniques, and in which doing so, one may not only find that orbits can easily get very complicated but also one may find that interesting yet also important statistical questions concerning measuring complexity of a dynamical system, in particular, the question of determining the topological entropy of a given system (which the reader need not worry about at all, since we won't be inspecting such a question), may become intractable. 

\begin{rem}
Loosely speaking, recall that \say{topological entropy} is a nonnegative statistic that gives a way of measuring the exponential growth rate of distinguishable orbits of a dynamical system as time progresses. You may see independently in more great details about topological entropy in the work of Adler \cite{Adl} and Bowen \cite{Ruf}.
\end{rem}Now before we proceed any further in our discussion, we wish to first take a very quick look at the following classical example of a polynomial morphism with a rational periodic point and rational periodic orbit: 

\begin{exmp} \label{ex 1.1}
Consider $\varphi_{2, -21/16}$ defined by $\varphi_{2, -21/16}(z) = z^2 -21/16$ for all $c,z \in\mathbb{Q}$. Then if we compute all iterations of $\varphi_{2, -21/16}$ on $z = 1/4$, we get that $\mathcal{O}^{+}(1/4)$ consists of points $\varphi^0(1/4) = 1/4$, $\varphi^1(1/4) = -5/4$, $\varphi^2(1/4) = 1/4$ and so on; and moreover the points in $\mathcal{O}^{+}(1/4)$ form a sequence of the following form:
\begin{center}
$\{1/4 \longrightarrow -5/4 \} \longrightarrow \{1/4 \longrightarrow -5/4\}\longrightarrow \{1/4 \longrightarrow -5/4\} \longrightarrow \cdots$ \end{center}So in this example, $z = 1/4$ is a periodic point of a map $\varphi_{2, -21/16}$ with exact period 2 and from the above definition, $\mathcal{O}^{+}(1/4)$ is a rational periodic orbit (i.e., $\mathcal{O}^{+}(1/4)$ is a periodic orbit consisting of rational points).
\end{exmp} 

In addition to the notion of a periodic point and a periodic orbit, we also have  in dynamical systems a more complicated but somewhat related notion of a preperiodic point and a preperiodic orbit. We call a point $\alpha\in {\mathbb{P}^N(K)}$ a \textit{preperiodic point of $\varphi$}, whenever $\varphi^{m+n}(\alpha) = \varphi^{m}(\alpha)$ for some integers $m\geq 0$ and $n\geq 1$. In this case, the smallest integers $m\geq 0$ and $n\geq 1$ such that $\varphi^{m+n}(\alpha) = \varphi^{m}(\alpha)$ happens, are called the \textit{preperiod} and \textit{eventual period of $\alpha$}, respectively. Again, we denote the set of preperiodic points of $\varphi$ by PrePer$(\varphi, {\mathbb{P}^N(K)})$. For any given preperiodic point $\alpha$ of $\varphi$, we then call the set of all iterates of $\varphi$ on $\alpha$, \textit{the preperiodic orbit of $\alpha$}.
Now observe that for $m=0$, we then have $\varphi^{n}(\alpha) = \alpha $ and so $\alpha$ is a periodic point of period $n$. Hence, every periodic point of $\varphi$ is also a preperiodic point of $\varphi$, that is, Per$(\varphi, {\mathbb{P}^N(K)}) \subseteq$ PrePer$(\varphi, {\mathbb{P}^N(K)})$; however, it need not be that the set PrePer$(\varphi, {\mathbb{P}^N(K)})\subseteq$ Per$(\varphi, {\mathbb{P}^N(K)})$ as illustrated by the following classical example: 
\begin{exmp}
Consider  $\varphi_{2, -29/16}$ defined by $\varphi_{2, -29/16}(z) = z^2 -29/16$ for all $c,z\in \mathbb{Q}$. Then if we compute all iterations of $\varphi_{2, -29/16}$ on $z = 3/4$, we get that $\mathcal{O}^{+}(3/4)$ consists of points $\varphi^0(3/4) = 3/4$, $\varphi^1(3/4) = -5/4$, $\varphi^2(3/4) = -1/4$ and so on; and moreover the points in $\mathcal{O}^{+}(3/4)$ form a sequence of the following form: 
\begin{center}
$3/4 \longrightarrow -5/4 \longrightarrow \{-1/4 \longrightarrow -7/4 \longrightarrow 5/4\} \longrightarrow \{-1/4 \longrightarrow -7/4 \longrightarrow 5/4\}\longrightarrow \cdots$    
\end{center}
So in this example, we conclude that $z = 3/4$ is a rational preperiodic point of $\varphi_{2, -29/16}$ with preperiod $m=2$ and eventual period $n=3$; and $\mathcal{O}^{+}(3/4)$ is a preperiodic orbit. Also, since preperiod $m= 2 \neq 0$, then $z = 3/4$ is not a periodic point. As another example of a preperiodic point, recall in Example \ref{ex 1.1} that $z = 1/4$ is a periodic point of $\varphi_{2, -21/16}$, and since every periodic point is a preperiodic point, then $z = 1/4$ is a $\mathbb{Q}$-preperiodic point of $\varphi_{2, -21/16}$. Other interesting examples of $\mathbb{Q}$-preperiodic points may be found in Poonen's work \cite{Poonen}.
\end{exmp}

In the year 1950, Northcott \cite{North} used the theory of height functions to show that not only is the set PrePer$(\varphi, {\mathbb{P}^N(K)})$ always finite, but also for a given morphism $\varphi$ the set PrePer$(\varphi, {\mathbb{P}^N(K)})$ can be computed effectively. Forty-five years later, in the year 1995, Morton and Silverman conjectured that PrePer$(\varphi, \mathbb{P}^N(K))$ can be bounded in terms of degree $d$ of $\varphi$, degree $D$ of $K$, and dimension $N$ of the space ${\mathbb{P}^N(K)}$. This celebrated conjecture is called the \textit{Uniform Boundedness Conjecture}; which we then restate here as the following conjecture:

\begin{conj} \label{silver-morton}[\cite{Morton}]
Fix integers $D \geq 1$, $N \geq 1$, and $d \geq 2$. There exists a constant $C'= C'(D, N, d)$ such that for all number fields $K/{\mathbb{Q}}$ of degree at most $D$, and all morphisms $\varphi: {\mathbb{P}^N}(K) \rightarrow {\mathbb{P}^N}(K)$ of degree $d$ defined over $K$, the total number of preperiodic points of a morphism $\varphi$ is at most $C'$, i.e., \#PrePer$(\varphi, \mathbb{P}^N(K)) \leq C'$.
\end{conj}
\noindent Note that a special case of Conjecture \ref{silver-morton} is when the degree $D$ of a number field $K$ is $D = 1$, dimension $N$ of a space $\mathbb{P}^N(K)$ is $N = 1$, and degree $d$ of a morphism $\varphi$ is $d = 2$. In this case, if $\varphi$ is a polynomial morphism, then it is a quadratic map defined over the field $\mathbb{Q}$. Moreover, in this very special case, in the year 1995, Flynn and Poonen and Schaefer conjectured that a quadratic map has no points $z\in\mathbb{Q}$ with exact period more than 3. This conjecture of Flynn-Poonen-Schaefer \cite{Flynn} (which has been resolved for cases $n = 4$, $5$ in \cite{mor, Flynn} respectively and conditionally for $n=6$ in \cite{Stoll} is, however, still open for all cases $n\geq 7$ and moreover, which also Hutz-Ingram \cite{Ingram} gave strong computational evidence supporting it) is restated here formally as the following conjecture. Note that in this same special case, rational points of exact period $n\in \{1, 2, 3\}$ were first found in the year 1994 by Russo-Walde \cite{Russo} and also found in the year 1995 by Poonen \cite{Poonen} using a different set of techniques. We restate here the anticipated conjecture of Flynn-Poonen-Schaefer as the following conjecture: 
 
\begin{conj} \label{conj:2.4.1}[\cite{Flynn}, Conjecture 2]
If $n \geq 4$, then there is no quadratic polynomial $\varphi_{2,c }(z) = z^2 + c\in \mathbb{Q}[z]$ with a rational point of exact period $n$.
\end{conj}
Now by assuming Conjecture \ref{conj:2.4.1} and also establishing interesting results on rational preperiodic points, in the year 1998, Poonen \cite{Poonen} then concluded that the total number of rational preperiodic points of any quadratic polynomial $\varphi_{2, c}(z)$ over $\mathbb{Q}$ is at most nine. We restate here formally Poonen's result as the following corollary:
\begin{cor}\label{cor2}[\cite{Poonen}, Corollary 1]
If Conjecture \ref{conj:2.4.1} holds, then $\#$PrePer$(\varphi_{2,c}, \mathbb{Q}) \leq 9$,  for all quadratic maps $\varphi_{2, c}$ defined by $\varphi_{2, c}(z) = z^2 + c$ for all points $c, z\in\mathbb{Q}$.
\end{cor}

Since Per$(\varphi, {\mathbb{P}^N(K)}) \subseteq$ PrePer$(\varphi, {\mathbb{P}^N(K)})$ and so if $\#$PrePer$(\varphi, \mathbb{P}^N(K))$ is bounded above, then so is $\#$Per$(\varphi, \mathbb{P}^N(K))$  also bounded above by the same upper bound. So now, we may extract out the following version \ref{per}; and this is because in Section \ref{sec2} and \ref{sec3} we study a dynamical setting in which $K$ is replaced with a ring $\mathbb{Z}$; hoping to understand (and not to prove) the possibility and validity of a periodic version of Conject.\ref{silver-morton}:

\begin{conj} \label{silver-morton 1}($(D, N) = (1,1)$-version of Conjecture \ref{silver-morton})\label{per}
Fix an integer $d \geq 2$. There exists a constant $C'= C'(d)$ such that for all morphisms $\varphi: {\mathbb{P}}^{1}(\mathbb{Q}) \rightarrow {\mathbb{P}}^{1}(\mathbb{Q})$ of degree $d$, the number $\#$Per$(\varphi, {\mathbb{P}}^{1}(\mathbb{Q})) \leq C'(d)$.
\end{conj}

\subsection*{History on the Connection Between the Size of Per$(\varphi_{d, c}, K)$ and the Coefficient $c$}

In the year 1994, Walde and Russo not only proved [\cite{Russo}, Corollary 4] that for a quadratic map $\varphi_{2,c}$ defined over $\mathbb{Q}$ with a periodic point, the denominator of a rational point $c$, denoted as den$(c)$, is a square but they also proved that den$(c)$ is even, whenever $\varphi_{2,c}$ admits a rational cycle of length $\ell \geq 3$. Moreover, Walde-Russo also proved [\cite{Russo}, Cor. 6, Thm 8 and Cor. 7] that the size \#Per$(\varphi_{2, c}, \mathbb{Q})\leq 2$, whenever den$(c)$ is an odd integer. 

Three years later, in the year 1997, Call-Goldstine \cite{Call} proved that the size of the set PrePer$(\varphi_{2,c},\mathbb{Q})$ of rational preperiodic points of a quadratic map $\varphi_{2, c}$ can be bounded above in terms of the number of distinct odd primes dividing the denominator of a point $c\in\mathbb{Q}$. We restate here formally this result of Call and Goldstine as the following theorem. Note that in this theorem, $GCD(a, e)$ refers to the greatest common divisor of $a$, $e \in \mathbb{Z}$:

\begin{thm}\label{2.3.1}[\cite{Call}, Theorem 6.9]
Let $e>0$ be an integer and let $s$ be the number of distinct odd prime factors of e. Define $\varepsilon  = 0$, $1$, $2$, if $4\nmid e$, if $4\mid e$ and $8 \nmid e$, if $8 \mid e$, respectively. Let $c = a/e^2$, where $a\in \mathbb{Z}$ and $GCD(a, e) = 1$. If $c \neq -2$, then the total number of $\mathbb{Q}$-preperiodic points of $\varphi_{2, c}$ is at most $2^{s + 2 + \varepsilon} + 1$. Moreover, a quadratic map $\varphi_{2, -2}$ has exactly six rational preperiodic points.
\end{thm}
\noindent
Now recall that the set Per$(\varphi_{2,c}, \mathbb{Q})$ is always a subset of PrePer$(\varphi_{2,c}, \mathbb{Q})$, and since we now know from Theorem \ref{2.3.1} that in the case of quadratic maps $\varphi_{2,c}$ defined over $\mathbb{Q}$, Call and Goldstine have shown that the maximum number of elements in PrePer$(\varphi_{2,c}, \mathbb{Q})$ is at most $2^{s + 2 + \varepsilon} + 1$ for any rational point $c\neq -2$. Hence, for any $c\neq -2$, the size of Per$(\varphi_{2,c}, \mathbb{Q})$ is equal to $2^{s + 2 + \varepsilon} + 1$ or is strictly less than $2^{s + 2 + \varepsilon} + 1$; and so for any $c\neq -2$, the size of Per$(\varphi_{2,c}, \mathbb{Q})$ is bounded above by a constant $2^{s + 2 + \varepsilon} + 1$ depending only on the number of distinct odd prime factors of den$(c)$ and on the quantity $\varepsilon$. On the other hand, for $c = -2$ (in which case den$(c) = 1$ and so has no distinct prime factors), Theorem \ref {2.3.1} tells us that \#PrePer$(\varphi_{2,-2}, \mathbb{Q})  = 6$ and so \#Per$(\varphi_{2,-2}, \mathbb{Q})  \leq 6$. 

Eight years later, after the work of Call-Goldstine, in the year 2005, Benedetto \cite{detto} conducted a detailed analysis of the filled Julia set and then improved substantially the bounds on the size of PrePre$(\varphi, K)$ which had been established by earlier work of Call-Goldstine and earlier works of several other people in the literature. In particular, by working only with polynomial maps $\varphi$ and only in dimension $N=1$, however, allowing arbitrary degree $d\geq 2$ and an arbitrary global field $K$ (i.e., $K$ is a number field or a function field defined over a finite field), Benedetto established the following remarkable result on the relationship between the size of the set PrePre$(\varphi, K)$ and the number of bad primes of $\varphi$ in $K$ (note that when the morphism $\varphi = \varphi_{d,c}$ for any point $c\in K$, then the \textit{bad primes} of $\varphi_{d,c}$ in $K$ are those primes that in fact divide the denominator of a $K$-point $c$):

\begin{thm}\label{main} [\cite{detto}, Main Theorem]
Let $K$ be a global field, $\varphi\in K[z]$ be a polynomial of degree $d\geq 2$ and $s$ be the number of bad primes of $\varphi$ in $K$. The number of preperiodic points of $\varphi$ in $\mathbb{P}^N(K)$ is at most $O(\text{s log s})$. 
\end{thm}

\noindent
Since from our earlier discussion we know that the set of periodic points of any given morphism $\varphi$ is always a subset of the set of preperiodic points of $\varphi$, so then, if $\varphi$ is a polynomial morphism of degree $d\geq 2$ defined over a global field $K$, then Thm \ref{main} also shows that the total number of periodic points of the underlying morphism $\varphi$ is in fact also $\ll \text{s log s}$, i.e., at most a positive constant times $\text{s log s}$. Notice that setting $\varphi(z) = \varphi_{d, c}(z) = z^d + c\in K[z]$ and $K = \mathbb{Q}$, then by Thm \ref{main} the size of Per$(\varphi_{d, c}, \mathbb{Q})$ is $\ll \text{s log s}$. More to this, since Thm \ref{main} applies to any $\varphi(z) = \varphi_{d, c}(z)$ of arbitrary degree $d\geq 2$ defined over $\mathbb{Q}$, it then follows that Thm \ref{main} can be applied to any $\varphi_{d, c}(z)$ of arbitrary prime or even degree $d\geq 2$; and as such one must then obtain the upper bound in Thm \ref{main} on the number of $\mathbb{Q}$-periodic points and hence on the number of $\mathbb{Z}$-periodic points. Hence, we again we see from Thm \ref{main} that the size of the set Per$(\varphi_{d, c}, \mathbb{Q})$ and thus of Per$(\varphi_{d, c}, \mathbb{Z})$, is bounded above using arithmetic information on the prime divisors of the coefficients of the polynomial $\varphi_{d, c}(z)$. 

Seven years after the work of Benedetto, in the year 2012, Narkiewicz's work \cite{Narkie} not only showed that any $\varphi_{d,c}$ defined over $\mathbb{Q}$ with odd degree $d\geq 3$ has no rational periodic points of exact period $n > 1$, but his also showed that the total number of $\mathbb{Q}$-preperiodic points is at most 4. We restate this result here as the following: 

\begin{thm} \label{theorem 3.2.1}\cite{Narkie}
For any integer $n > 1$ and any odd integer $d\geq 3$, there is no $c\in \mathbb{Q}$ such that $\varphi_{d,c}$ defined by $\varphi_{d, c}(z)$ for all $c,z \in \mathbb{Q}$ has rational periodic points of exact period $n$. Moreover, $\#PrePer(\varphi_{d, c}, \mathbb{Q}) \leq 4$. 
\end{thm} 

\begin{rem}
The first part of Thm \ref{theorem 3.2.1} is proved by observing that for each odd degree $d\geq 3$, the polynomial $\varphi_{d, c}(z)$ is non-decreasing on $\mathbb{Q}$ and so by elementary mathematical analysis one then expects the forward orbit $\mathcal{O}_{\varphi_{d,c}} (z)$ of each rational point $z$ to form a non-decreasing sequence of iterations. Hence, it is immediately evident that the polynomial $\varphi_{d, c}(z)$ can only have rational points of exact period $n = 1$ (with no preperiod). The upper bound 4 is obtained by counting the number of rational roots of $\varphi_{d, c}(z)-z$. Notice Thm \ref{theorem 3.2.1} shows that the number of rational points $z$ that satisfy the equation $z^d - z + c = 0$ is bounded above by 4; which then also means that the number of such rational points and hence such $\mathbb{Z}$-points is equal to 4 or strictly less than 4. 
\end{rem}

Three years after \cite{Narkie}, in 2015, Hutz \cite{Hutz} developed an algorithm determining effectively all $\mathbb{Q}$-preperiodic points of a morphism defined over a given number field $K$; from which he then made the following conjecture:

\begin{conj} \label{conjecture 3.2.1}[\cite{Hutz}, Conjecture 1a]
For any integer $n > 2$, there is no even degree $d > 2$ and no point $c \in \mathbb{Q}$ such that the polynomial map $\varphi_{d, c}$ has rational points of exact period $n$.
Moreover, \#PrePer$(\varphi_{d, c}, \mathbb{Q}) \leq 4$. 
\end{conj}

\begin{rem} \label{nt}
If Conjecture \ref{conjecture 3.2.1} held, it would then follow that the total number of integral fixed points of any $\varphi_{d,c}(x)$ of even degree $d>2$ is equal to $4$ or $<4$. Moreover, since the monic polynomial $\varphi_{d,c}(x)\in \mathbb{Z}[x]$ has good reduction modulo $p$, it then follows that the total number of integral fixed points of any $\varphi_{d,c}(x)$ modulo $p$ is also $4$ or $<4$. But of course now the issue is that we unfortunately don't know (as to the author's knowledge) whether Conjecture \ref{conjecture 3.2.1} holds or not, let alone whether $4$ is the correct upper bound on the total number of rational (and hence integral) fixed points of any $\varphi_{d,c}$ of even degree $d>2$. On whether any progress has been made on Conjecture \ref{conjecture 3.2.1}, recently Panraksa \cite{par2} proved among many other results that $\varphi_{4,c}(z)\in\mathbb{Q}[z]$ has $\mathbb{Q}$-points of exact period $n = 2$. Moreover, he also proved that $\varphi_{d,c}(z)\in\mathbb{Q}[z]$ has no $\mathbb{Q}$-points of exact period $n = 2$ for any $c \in \mathbb{Q}$ with $c \neq -1$ and $d = 6$, $2k$ with $3 \mid 2k-1$. The interested reader may find these mentioned results of Panraksa in his unconditional Thms. 2.1, 2.4 and Thm. 1.7 conditioned on the abc-conjecture in \cite{par2}.
\end{rem}

Twenty-eight years later, after the work of Walde-Russo, in the year 2022, Eliahou and Fares proved [\cite{Shalom2}, Theorem 2.12] that the denominator of a rational point $-c$, denoted as den$(-c)$ is divisible by 16, whenever a quadratic map $\varphi_{2,-c}$ defined by $\varphi_{2, -c}(z) = z^2 - c$ for all points $c, z\in \mathbb{Q}$ admits a rational cycle of length $\ell \geq 3$. Moreover, Eliahou and Fares also proved [\cite{Shalom2}, Proposition 2.8] that the size \#Per$(\varphi_{2, -c}, \mathbb{Q})\leq 2$, whenever den$(-c)$ is an odd integer. Motivated by also the work of Call-Goldstine, Eliahou-Fares \cite{Shalom2} also proved that the size of the set Per$(\varphi_{2, -c}, \mathbb{Q})$ can be bounded above simply by using information on den$(-c)$, namely, information in terms of the number of distinct primes dividing den$(-c)$. Moreover, the authors \cite{Shalom1} also showed that the upper bound is 4, whenever $c\in \mathbb{Q^*} = \mathbb{Q}\setminus\{0\}$. We restate here their results as the following:

\begin{cor}\label{sha}[\cite{Shalom2, Shalom1}, Corollary 3.11 and Corollary 4.4, resp.]
Let $c\in \mathbb{Q}$ such that den$(c) = d^2$ with $d\in 4 \mathbb{N}$. Let $s$ be the number of distinct primes dividing $d$. Then, the total number of $\mathbb{Q}$-periodic points of $\varphi_{2, -c}$ is at most $2^s + 2$. Moreover, for $c\in \mathbb{Q^*}$ such that the den$(c)$ is a power of a prime number. Then, $\#$Per$(\varphi_{2, c}, \mathbb{Q}) \leq 4$.
\end{cor}

\begin{rem}
Since every point in a given $\mathbb{Q}$-periodic cycle of $\varphi_{2,c}$ is again a $\mathbb{Q}$-periodic point of $\varphi_{2,c}$, then it seems somewhat reasonable to infer from the foregoing history that what Walde-Russo initially proved and later improved by Eliahou-Fares is, precisely the claim that den$(c)$ is divisible by 2, whenever $\#$Per$(\varphi_{2, c}, \mathbb{Q})\geq 3$ or $\#$Per$(\varphi_{2, -c}, \mathbb{Q})\geq 3$. On the other hand, what both Walde-Russo and Eliahou-Fares also proved is the claim that the number $\#$Per$(\varphi_{2, c}, \mathbb{Q})\leq 2$ or $\#$Per$(\varphi_{2, -c}, \mathbb{Q})\leq 2$, respectively, whenever den$(c)$ is an odd integer. So now, it's highly plausible to assert that what the authors in \cite{Russo, Call, Shalom1, Shalom2} have studied is this very surprising relationship between the total number of elements in Per$(\varphi_{2, c}, \mathbb{Q})$ or Per$(\varphi_{2, -c},\mathbb{Q})$ and denominator of a generic coefficient $c$ or $-c$ of $\varphi_{2, c}(z)$ or $\varphi_{2, -c}(z)$, respectively. For $d\geq 2$, we've seen that what the authors in \cite{detto, Narkie, Hutz, par2} studied, is this same relationship between size of Per$(\varphi_{d, c}, \mathbb{Q})$ and the coefficients of $\varphi_{d, c}(z)$.
\end{rem}

\noindent The purpose of this article is once again to inspect further the above connection, independently in the case of polynomial maps $\varphi_{p, c}$ of odd prime degree $p$ defined over $\mathbb{Z}$ for any given prime integer $p\geq 3$ and in the case of polynomial maps $\varphi_{p-1, c}$ of even degree $p-1$ defined over $\mathbb{Z}$ for any given prime integer $p\geq 5$; and doing so from a spirit that is truly inspired and guided by some of the many striking developments in arithmetic statistics.

\section{On the Number of Integral Fixed Points of any Family of Polynomial Maps $\varphi_{p,c}$}\label{sec2}

In this section, we use an elementary number-theoretic fact combined with a good reduction fact applied on Theorem \ref{theorem 3.2.1} to then count the number of distinct fixed points of any $\varphi_{p,c}$ modulo $p$, where $p\geq 3$ is any prime. To do so, we let $c\in \mathbb{Z}$ be any integer and $p\geq 3$ be any prime, and then define fixed point-counting function 
\begin{equation}\label{N_{c}}
N_{c}(p) := \# \{ z\in \mathbb{Z} / p\mathbb{Z} : \varphi_{p,c}(z) - z \equiv 0 \ \text{(mod $p$)}\}.
\end{equation}\noindent We then first prove the following theorem on number of distinct fixed points of any cubic map $\varphi_{3, c}$ modulo $3$:
\begin{thm} \label{2.1}
Given any family of cubic maps $\varphi_{3, c}$ defined by $\varphi_{3, c}(z) = z^3 + c$ for all $c, z\in\mathbb{Z}$, and let $N_{c}(3)$ be defined as in \textnormal{(\ref{N_{c}})}. Then $N_{c}(3) = 3$ for every coefficient $c = 3t$; otherwise $N_{c}(3) = 0$ for every coefficient $c \neq  3t$.
\end{thm}
\begin{proof}
Let $f(z)= \varphi_{3, c}(z)-z = z^3 - z + c$ and note that for every coefficient $c = 3t$, we then write $f(z)= z^3 - z + 3t$. Now reducing $f(z)$ modulo $3$, it then follows $f(z)\equiv z^3 - z$ (mod $3$). But then for every integer $z\in \mathbb{Z}$, it then follows by Fermat's Little Theorem that $z^3-z\equiv 0$ (mod $3$); which also means $z(z^2-1)\equiv 0$ (mod $3$). It then follows that either $z\equiv 0$ (mod $3$) or $z^2 - 1\equiv 0$ (mod $3$), since $\mathbb{Z}\slash 3\mathbb{Z}$ is an integral domain. Moreover, by a standard fact in elementary number theory, we note that $z^2 -1 \equiv 0$ (mod $3$) implies $z\equiv \pm 1$ (mod $3$). This then means $f(z)=0$ has exactly three unique solutions $z\equiv -1, 0, 1$ (mod $3$); and so every integer solution of $f(z)=0$ is of the form $z_{1} = 3r-1$, $z_{2} = 3s$ and $z_{3} = 3u+1$ (mod $3$) for some $r, s, u\in \mathbb{Z}$. But now we then conclude $\#\{ z\in \mathbb{Z} / 3\mathbb{Z} : \varphi_{3,c}(z) - z \equiv 0 \text{ (mod 3)}\} = 3$ and so $N_{c}(3) = 3$. To see $N_{c}(3) = 0$ for every coefficient $c \neq  3t$, we first note that since $c\neq 3t$, then this also means $c\not \equiv 0$ (mod $3$). But now using the fact that $z^3 = z$ for every element $z\in \mathbb{Z}\slash 3\mathbb{Z}$, we then note that $z^3 - z + c\not \equiv 0$ (mod $3$) for every point $z\in \mathbb{Z} / 3\mathbb{Z}$; and so $f(z)\not \equiv 0$ (mod $3$) for every point $z\in \mathbb{Z} / 3\mathbb{Z}$. This then means $f(x)=\varphi_{3,c}(x)-x$ has no roots in $\mathbb{Z} / 3\mathbb{Z}$ for every coefficient $c$ not divisible by $3$; and thus  we then conclude $N_{c}(3) = 0$. This then completes the whole proof, as desired.
\end{proof} 
We now wish to generalize Theorem \ref{2.1} to any $\varphi_{p, c}$ for any given prime $p\geq 3$. More precisely, assuming Theorem \ref{theorem 3.2.1}, we then prove that the number of distinct fixed points of any $\varphi_{p, c}$ modulo $p$ is also either $3$ or $0$:

\begin{thm} \label{2.2}
Let $p\geq 3$ be any fixed prime integer, and assume Theorem \ref{theorem 3.2.1}. Consider any family of polynomial maps $\varphi_{p, c}$ defined by $\varphi_{p, c}(z) = z^p + c$ for all points $c, z\in\mathbb{Z}$, and let $N_{c}(p)$ be the number defined as in \textnormal{(\ref{N_{c}})}. Then $N_{c}(p) = 3$ for every coefficient $c = pt$; otherwise the number $N_{c}(p) = 0$ for any coefficient  $c \neq  pt$.
\end{thm}
\begin{proof}
By applying a somewhat similar argument as in the Proof of Theorem \ref{2.1}, we then obtain the count as desired. That is, let $f(z)= \varphi_{p,c}(z)-z=z^p - z + c$ and for every coefficient $c = pt$, we then write $f(z)= z^p - z + pt$. Now reducing $f(z)$ modulo $p$, it then follows $f(z)\equiv z^p - z$ (mod $p$). But then recall by Fermat's Little Theorem that $z^p-z\equiv 0$ (mod $p$) for every $z\in \mathbb{Z}$; equivalently $z(z^{p-1}-1)\equiv 0$ (mod $p$) for every $z\in \mathbb{Z}$. This then means $z\equiv 0$ (mod $p$) or $z^{p-1} - 1\equiv 0$ (mod $p$), since $\mathbb{Z} / p\mathbb{Z}$ is an integral domain. Moreover, we note by a standard fact that $z^{p-1} -1 \equiv 0$ (mod $p$) holds for every $z\in \{1, 2, \cdots, p-1 \}$; and so $f(z)\equiv 0$ (mod $p$) for every $z\in \{1, 2, \cdots, p-1 \}$. But now recall from Theorem \ref{theorem 3.2.1} that $\#\{ z\in \mathbb{Z} : \varphi_{p,c}(z) - z = 0\}\leq 4$, and since by [\cite{Silverman}, Corollary 2.20] we know that good reduction modulo $p$ preserves periodicity of points, then together with $z\equiv 0$ (mod $p$), we note that at most three additional values of $z\in \{1, 2, \cdots, p-1 \}$ are desired for $f(x)\equiv 0$ (mod $p$). But now because of the first part of Theorem \ref{2.1}, we note that for every $p\geq 3$, then $f(z)$ modulo $p$ has three roots in $\mathbb{Z} / p\mathbb{Z}$. But now we then conclude $\#\{ z\in \mathbb{Z} / p\mathbb{Z} : \varphi_{p,c}(z) - z \equiv 0 \text{ (mod $p$)}\} = 3$ and so $N_{c}(p) = 3$. To see $N_{c}(p) = 0$ for any coefficient  $c \neq  pt$, we note that since $c\neq pt$, then $c\not \equiv 0$ (mod $p$). But now using (as a fact) that $z^p = z$ for every $z\in \mathbb{Z}\slash p \mathbb{Z}$, we then note $z^p - z + c\not \equiv 0$ (mod $p$) for every $z\in \mathbb{Z} / p\mathbb{Z}$; and so $f(z)\not \equiv 0$ (mod $p$) for every $z\in \mathbb{Z} / p\mathbb{Z}$. This then means $f(x)=\varphi_{p,c}(x)-x$ has no roots in $\mathbb{Z} / p\mathbb{Z}$ for every coefficient $c\neq pt$; and so we conclude $N_{c}(p) = 0$. This then completes the whole proof, as required.
\end{proof}

\begin{rem}\label{rem2.3}
With now Theorem \ref{2.2}, we may then to each distinct integral fixed point of $\varphi_{p,c}$ associate an integral fixed orbit. In doing so, we then obtain a dynamical translation of Theorem \ref{2.2}, namely, that the number of distinct integral fixed orbits that any $\varphi_{p,c}$ has when iterated on $\mathbb{Z} / p\mathbb{Z}$ is $3$ or $0$. Notice that in both coefficient cases $c\equiv 0$ (mod $p)$ and $c\not \equiv 0$ (mod $p)$ that we considered in Theorem \ref{2.2}, it may also follow that the expected total number of distinct integral fixed points in the whole family of maps $\varphi_{p,c}$ modulo $p$ is $3 + 0 =3$.
\end{rem}
\section{The Number of Integral Fixed Points of any Family of Polynomial Maps $\varphi_{p-1,c}$}\label{sec3}

Unlike in Section \ref{sec2} in which we are assumed a theorem, we in this section also wish to count unconditionally the number of distinct integral fixed points of any $\varphi_{p-1,c}$ modulo $p$, where $p\geq 5$ is any given prime. To that end, we again let $c\in \mathbb{Z}$ be any integer and $p\geq 5$ be any prime, and then define fixed point-counting function 
\begin{equation}\label{M_{c}}
M_{c}(p) := \#\{ z\in \mathbb{Z} / p\mathbb{Z} : \varphi_{p-1,c}(z) - z \equiv 0 \ \text{(mod $p$)}\}.
\end{equation}\noindent As before, we then first prove the following theorem on number of distinct fixed points of any $\varphi_{4, c}$ modulo $5$:
\begin{thm} \label{6.0.1}
Let $\varphi_{4, c}$ be defined by $\varphi_{4, c}(z) = z^4 + c$ for all $c, z\in\mathbb{Z}$, and $M_{c}(5)$ be defined as in \textnormal{(\ref{M_{c}})}. Then $M_{c}(5) = 1$ or $2$ for any coefficient $c\equiv 1 \ (mod \ 5)$ or $c = 5t$, resp.; otherwise $M_{c}(5) = 0$ for any $c\equiv -1\ (mod \ 5)$. 
\end{thm}
\begin{proof}
Let $g(z)= \varphi_{4,c}(z)-z = z^4 - z + c$ and note that for every coefficient $c = 5t$, reducing $g(z)$ modulo $5$, we then obtain $g(z)\equiv z^4 - z$ (mod $5$); and so the reduced polynomial $g(z)$ modulo $5$ is now a polynomial defined over a finite  field $\mathbb{Z}\slash5\mathbb{Z}$ of $5$ distinct elements. So now, it is a well-known fact about polynomials over finite fields that the quartic monic polynomial $h(x)=x^4 -1$ has $4$ distinct nonzero roots in $\mathbb{Z}\slash5\mathbb{Z}$; and so $z^4 = 1$ for every nonzero element $z\in \mathbb{Z}\slash5\mathbb{Z}$. But then $g(z)\equiv 1 - z$ (mod $5$) for every nonzero point $z\in \mathbb{Z}\slash5\mathbb{Z}$; from which it then follows that the polynomial $g(z)$ has a nonzero root in $\mathbb{Z}\slash5\mathbb{Z}$, namely, $z\equiv 1$ (mod $5$). Moreover, since $z$ is also a linear factor of $g(z)\equiv z(z^3 - 1)$ (mod $5$), it then follows $z\equiv 0$ (mod $5$) is also a root of $g(z)$ modulo $5$. But now we then conclude $\#\{ z\in \mathbb{Z} / 5\mathbb{Z}: \varphi_{4,c}(z) - z \equiv 0 \text{ (mod 5)}\} = 2$ for every coefficient $c\equiv 0$ (mod $5$) and so the number $M_{c}(5) = 2$. To see that the number $M_{c}(5) = 1$ for every coefficient $c\equiv 1$ (mod $5$), we first note that since $c\equiv 1$ (mod $5$) and since also $z^4 = 1$ for every nonzero element $z\in \mathbb{Z}\slash5\mathbb{Z}$, then reducing $g(z)= \varphi_{4,c}(z)-z = z^4 - z + c$ modulo $5$, it then follows $g(z)\equiv 2 - z$ (mod $5$) and so $g(z)$ has a root in $\mathbb{Z}\slash5\mathbb{Z}$, namely, $z\equiv 2$ (mod $5$); and so we conclude $M_{c}(5) = 1$. To see that the number $M_{c}(5) = 0$ for every coefficient $c\equiv -1$ (mod $5$), we note that since $c \equiv -1$ (mod $5$) and since also (as a fact) $z^4 = 1$ for every nonzero $z\in \mathbb{Z}\slash5\mathbb{Z}$, it then follows $z^4 - z + c\equiv -z$ (mod $5$); and so $g(z) \equiv -z$ (mod $5$). But now notice that $z\equiv 0$ (mod $5$) is a root of $g(z)$ modulo $5$ for every coefficient $c\equiv -1$ (mod $5$) and for every coefficient $c\equiv 0$ (mod $5$) as seen from the first part; which then clearly is impossible, since $-1\not \equiv 0$ (mod $5$). This then means that $g(x)=\varphi_{4,c}(x)-x$ has no roots in $\mathbb{Z} \slash 5\mathbb{Z}$ for every coefficient $c\equiv -1$ (mod $5$); and so we then conclude $M_{c}(5) = 0$, as desired.
\end{proof} 
We now wish to generalize Theorem \ref{6.0.1} to any polynomial map $\varphi_{p-1, c}$ for any given prime $p\geq 5$. Specifically, we prove that the number of distinct fixed points of every $\varphi_{p-1, c}$ modulo $p$ is equal to $1$ or $2$ or $0$:

\begin{thm} \label{6.0.2}
Let $p\geq 5$ be any fixed prime integer, and $\varphi_{p-1, c}$ be any polynomial map defined by $\varphi_{p-1, c}(z) = z^{p-1}+c$ for all $c, z\in\mathbb{Z}$. Let $M_{c}(p)$ be the number defined as in \textnormal{(\ref{M_{c}})}. Then we have $M_{c}(p) = 1$ or $2$ for every coefficient $c\equiv 1 \ (mod \ p)$ or $c = pt$, respectively; otherwise we have $M_{c}(p) = 0$ for every point $c\equiv -1\ (mod \ p)$.
\end{thm}
\begin{proof}
By applying a similar argument as in the Proof of Theorem \ref{6.0.1}, we then obtain the count as desired. That is, let $g(z)= \varphi_{p-1,c}(z)-z=z^{p-1} - z + c$ and for every coefficient $c = pt$, reducing $g(z)$ modulo $p$, it then follows $g(z)\equiv z^{p-1} - z$ (mod $p$); and so the reduced polynomial $g(z)$ modulo $p$ is now a polynomial defined over a finite field $\mathbb{Z}\slash p\mathbb{Z}$ of order $p$. So now, as before we know (as a fact) that the monic polynomial $h(x)=x^{p-1} -1$ has $p-1$ distinct nonzero roots in $\mathbb{Z}\slash p\mathbb{Z}$; and so $z^{p-1} = 1$ for every nonzero $z\in \mathbb{Z}\slash p\mathbb{Z}$. But then the reduced polynomial $g(z)\equiv 1 - z$ (mod $p$) for every nonzero point $z\in \mathbb{Z}\slash p\mathbb{Z}$; and so $g(z)$ modulo $p$ has a nonzero root in $\mathbb{Z}\slash p\mathbb{Z}$, namely, $z\equiv 1$ (mod $p$). Moreover, since $z$ is a linear factor of the reduced polynomial $g(z)\equiv z(z^{p-2} - 1)$ (mod $p$), it then follows $z\equiv 0$ (mod $p$) is also a root of the polynomial $g(z)$ modulo $p$. But now we then conclude $\#\{ z\in \mathbb{Z} / p\mathbb{Z} : \varphi_{p-1,c}(z) - z \equiv 0 \text{ (mod $p$)}\} = 2$ and so the number $M_{c}(p) = 2$. To see that the number $M_{c}(p) = 1$ for every coefficient $c\equiv 1$ (mod $p$), we note that since the coefficient $c\equiv 1$ (mod $p$) and since also $z^{p-1} = 1$ for every nonzero element $z\in \mathbb{Z}\slash p\mathbb{Z}$, then reducing $g(z)= \varphi_{p-1,c}(z)-z = z^{p-1} - z + c$ modulo $p$, it then follows $g(z)\equiv 2 - z$ (mod $p$); and so $g(z)$ modulo $p$ has a root in $\mathbb{Z}\slash p\mathbb{Z}$, namely, $z\equiv 2$ (mod $p$); and so we then conclude $M_{c}(p) = 1$. Finally, to see that the number $M_{c}(p) = 0$ for every coefficient $c\equiv -1$ (mod $p$), we as before note that since the coefficient $c \equiv -1$ (mod $p$) and since also (as a fact) $z^{p-1} = 1$ for every nonzero element $z\in \mathbb{Z}\slash p\mathbb{Z}$, it then follows $z^{p-1} - z + c\equiv -z$ (mod $p$); and so $g(z) \equiv -z$ (mod $p$). But now notice $z\equiv 0$ (mod $p$) is a root of $g(z)$ modulo $p$ for every coefficient $c\equiv -1$ (mod $p$) and every coefficient $c\equiv 0$ (mod $p$) as seen from the first part; which then is impossible, since $-1\not \equiv 0$ (mod $p$). It then follows $g(x)=\varphi_{p-1,c}(x)-x$ has no roots in $\mathbb{Z} \slash p\mathbb{Z}$ for every coefficient $c\equiv -1$ (mod $p$); and so we then conclude $M_{c}(p) = 0$, as desired.
\end{proof}

\begin{rem}
With now Theorem \ref{6.0.2}, we may then to each distinct integral fixed point of $\varphi_{p-1,c}$ associate an integral fixed orbit. In doing so, we then also obtain a dynamical translation of Theorem \ref{6.0.2}, namely, that the number of distinct integral fixed orbits that $\varphi_{p-1,c}$ has when iterated on the space $\mathbb{Z} / p\mathbb{Z}$ is equal to $1$ or $2$ or $0$. Again, observe that in all of the possible coefficient cases $c\equiv \pm 1, 0$ (mod $p)$ considered in Theorem \ref{6.0.2}, it may then also follow that the expected total number of distinct integral fixed points in the whole family of maps $\varphi_{p-1,c}$ modulo $p$ is equal to $1 + 2 + 0 =3$; which somewhat surprising coincides not only with the second prediction on the upper bound observed in Remark \ref{nt} on Conj.\ref{conjecture 3.2.1}, but also coincides with the expected total number three of distinct integral fixed points in the whole family of maps $\varphi_{p,c}$ modulo $p$ noted in Remark \ref{rem2.3}. 
\end{rem}

\section{On the Average Number of Fixed Points of any Polynomial Map $\varphi_{p,c}$ and $\varphi_{p-1,c}$}\label{sec4}

In this section, we wish to inspect independently the behavior of fixed point-counting functions $N_{c}(p)$ and $M_{c}(p)$ as $c$ tends to infinity. First, we wish to determine: \say{\textit{What is the average value of the function $N_{c}(p)$ as $c \to \infty$?}} The following corollary shows that the average value of the function $N_{c}(p)$ exits and is equal to $3$ or $0$ as $c\to \infty$:
\begin{cor}\label{4.1}
Let $p\geq 3$ be any prime integer. Then the average value of fixed point-counting function $N_{c}(p)$ exists and is equal to $3$ or $0$ as $c \to \infty$. More precisely, we have
\begin{myitemize}
    \item[\textnormal{(a)}] \textnormal{Avg} $N_{c = pt}(p) := \lim\limits_{c\to\infty} \Large{\frac{\sum\limits_{3\leq p\leq c, \ p\mid c}N_{c}(p)}{\Large{\sum\limits_{3\leq p\leq c, \ p\mid c}1}}} = 3.$  
    
    \item[\textnormal{(b)}] \textnormal{Avg} $N_{c\neq pt}(p):= \lim\limits_{c \to\infty} \Large{\frac{\sum\limits_{3\leq p\leq c, \ p\nmid c}N_{c}(p)}{\Large{\sum\limits_{3\leq p\leq c, \ p\nmid c}1}}} =  0$.    
\end{myitemize}

\end{cor}
\begin{proof}
Since from Theorem \ref{2.2} we know that the number $N_{c}(p) = 3$ for any prime $p\mid c$, we then obtain that $\lim\limits_{c\to\infty} \Large{\frac{\sum\limits_{3\leq p\leq c, \ p\mid c}N_{c}(p)}{\Large{\sum\limits_{3\leq p\leq c, \ p\mid c}1}}} = 3\lim\limits_{c\to\infty} \Large{\frac{\sum\limits_{3\leq p\leq c, \ p\mid c}1}{\Large{\sum\limits_{3\leq p\leq c, \ p\mid c}1}}} = 3$. Hence, the average value of $N_{c}(p)$, namely, Avg $N_{c = pt}(p)$ is equal to 3, which shows $\textnormal{(a)}$. To see $\textnormal{(b)}$, we recall from Theorem \ref{2.2} that $N_{c}(p) = 0$ for any prime $p\nmid c$; and so we then obtain $\lim\limits_{c\to\infty} \Large{\frac{\sum\limits_{3\leq p\leq c, \ p\nmid c}N_{c}(p)}{\Large{\sum\limits_{3\leq p\leq c, \ p\nmid c}1}}} = 0$; and so Avg $N_{c\neq pt}(p) = 0$. This then completes the whole proof, as desired.
\end{proof}
\begin{rem} \label{re5}
From arithmetic statistics to arithmetic dynamics, we note that Corollary \ref{4.1} shows that any $\varphi_{p,c}$ iterated on the space $\mathbb{Z} / p\mathbb{Z}$ has on average three or zero distinct integral fixed orbits as $c$ tends to infinity.
\end{rem}

Similarly, we also wish to determine: \say{\textit{What is the average value of the function $M_{c}(p)$ as $c \to \infty$?}} The following corollary shows that the average value of the function $M_{c}(p)$ exits and equal to $1$ or $2$ or $0$ as $c\to \infty$:
\begin{cor}\label{cor5}
Let $p\geq 5$ be any prime integer. Then the average value of fixed point-counting function $M_{c}(p)$ exists and is equal to $1$ or $2$ or $0$ as $c\to\infty$. Specifically, we have 
\begin{myitemize}
    \item[\textnormal{(a)}] \textnormal{Avg} $M_{c-1 = pt}(p) := \lim\limits_{c\to\infty} \Large{\frac{\sum\limits_{5\leq p\leq (c-1), \ p\mid (c-1)}M_{c}(p)}{\Large{\sum\limits_{5\leq p\leq (c-1), \ p\mid (c-1)}1}}} = 1.$ 

    \item[\textnormal{(b)}] \textnormal{Avg} $M_{c= pt}(p) := \lim\limits_{c\to\infty} \Large{\frac{\sum\limits_{5\leq p\leq c, \ p\mid c}M_{c}(p)}{\Large{\sum\limits_{5\leq p\leq c, \ p\mid c}1}}} = 2.$
    
    \item[\textnormal{(c)}] \textnormal{Avg} $M_{c+1= pt}(p):= \lim\limits_{c \to\infty} \Large{\frac{\sum\limits_{5\leq p\leq (c+1), \ p\mid (c+1)}M_{c}(p)}{\Large{\sum\limits_{5\leq p\leq (c+1), \ p\mid (c+1)}1}}} =  0$.    
\end{myitemize}

\end{cor}
\begin{proof}
By applying a similar argument as in the Proof of Corollary \ref{4.1}, we then obtain the limits, as desired.
\end{proof} 
\begin{rem} \label{re5}
As before, from arithmetic statistics to arithmetic dynamics, we also note that Corollary \ref{cor5} demonstrates that any $\varphi_{p-1,c}$ iterated on $\mathbb{Z} / p\mathbb{Z}$ has on average one or two or zero integral fixed orbits as $c\to \infty$.
\end{rem}
\section{On the Density of Monic Integer Polynomials $\varphi_{p,c}(x)$ with the number $N_{c}(p) = 3$}\label{sec5}
In this section, we wish to ask and answer: \say{\textit{For $p\geq 3$ a prime integer, what is the density of monic integer polynomials $\varphi_{p, c}(x) = x^p + c$ with exactly three integral fixed points modulo $p$?}} The following corollary shows that there are very few integer polynomials $\varphi_{p,c}(x) = x^p + c$ with exactly three integral fixed points modulo $p$:
\begin{cor}\label{Thm 4.1}
Let $p\geq 3$ be a prime integer. Then the density of monic polynomials $\varphi_{p,c}(x) = x^p + c\in \mathbb{Z}[x]$ with the number $N_{c}(p) = 3$ exists and is equal to $0 \%$ as $c\to \infty$. More precisely, we have 
\begin{center}
    $\lim\limits_{c\to\infty} \Large{\frac{\# \{\varphi_{p,c}(x)\in \mathbb{Z}[x] \ : \ 3\leq p\leq c \ \text{and} \ N_{c}(p) \ = \ 3\}}{\Large{\# \{\varphi_{p,c}(x) \in \mathbb{Z}[x] \ : \ 3\leq p\leq c \}}}} = \ 0.$
\end{center}
\end{cor}
\begin{proof}
Since the defining condition $N_{c}(p) = 3$ is as we proved in Theorem \ref{2.2}, determined whenever the coefficient $c$ is divisible by any prime $p\geq 3$, we may then count $\# \{\varphi_{p,c}(x) \in \mathbb{Z}[x] : 3\leq p\leq c \ \text{and} \ N_{c}(p) \ = \ 3\}$ by counting the number $\# \{\varphi_{p,c}(x)\in \mathbb{Z}[x] : 3\leq p\leq c \ \text{and} \ p\mid c \ \text{for \ any \ fixed} \ c \}$. In that case, we then write 
\begin{center}
$\Large{\frac{\# \{\varphi_{p,c}(x) \in \mathbb{Z}[x] \ : \ 3\leq p\leq c \ \text{and} \ N_{c}(p) \ = \ 3\}}{\Large{\# \{\varphi_{p,c}(x) \in \mathbb{Z}[x] \ : \ 3\leq p\leq c \}}}} = \Large{\frac{\# \{\varphi_{p,c}(x)\in \mathbb{Z}[x] \ : \ 3\leq p\leq c \ \text{and} \ p\mid c \ \text{for any fixed} \ c \}}{\Large{\# \{\varphi_{p,c}(x) \in \mathbb{Z}[x] \ : \ 3\leq p\leq c \}}}}$. 
\end{center}\indent Moreover, for any fixed integer $c\geq 3$, the numerator of the foregoing quotient may be rewritten to then obtain
\begin{center}
$\# \{\varphi_{p,c}(x) \in \mathbb{Z}[x] : 3\leq p\leq c \ \text{and} \ p\mid c \} = \# \{p : 3\leq p\leq c \text{ and } p\mid c \} = \sum_{3\leq p\leq c, \ p\mid c}1 = \omega (c)$, 
\end{center}where $\omega(n)$ is by definition the number of distinct prime factors of $n$. Writing $\# \{\varphi_{p,c}(x) \in \mathbb{Z}[x]  : 3\leq p\leq c \} = \sum_{3\leq p\leq c} 1 = \pi(c)$, where $\pi(n)$ is by definition the number of primes at most $n$, we then note that the quotient 
\begin{center}
$\Large{\frac{\# \{\varphi_{p,c}(x)\in \mathbb{Z}[x] \ : \ 3\leq p\leq c \ \text{and} \ p\mid c \ \text{for any fixed} \ c \}}{\Large{\# \{\varphi_{p,c}(x)\in \mathbb{Z}[x] \ : \ 3\leq p\leq c \}}}} = \frac{\omega(c)}{\pi(c)}$.
\end{center}So now, recall (from a well-known fact) that for any $c\in \mathbb{Z}_{\geq 3}$, we have $2^{\omega(c)}\leq \sigma (n) \leq 2^{\Omega(c)}$, where $\sigma(n)$ is by definition the divisor function and $\Omega(n)$ is by definition the total number of prime factors of $n$, with respect to their multiplicity. Note that taking logarithms, we then obtain $\omega(c)\leq \frac{\text{log} \ \sigma(c)}{\text{log} \ 2}$; and so $\frac{\omega(c)}{\pi(c)} \leq \frac{\text{log} \ \sigma(c)}{\text{log} \ 2 \cdot \pi(c)}$. Moreover, for every $\epsilon >0$, it is well-known that $\sigma(c) = o(c^{\epsilon})$; and so log $\sigma(c) =$ log $o(c^{\epsilon})$ and so have $\frac{\omega(c)}{\pi(c)} \leq \frac{\text{log} \ o(c^{\epsilon})}{\text{log} \ 2 \cdot \pi(c)}$. Now for every fixed $\epsilon>0$, we then note $\lim\limits_{c\to\infty} \frac{\text{log} \ o(c^{\epsilon})}{\text{log} \ 2 \cdot \pi(c)} = 0$ and so $\lim\limits_{c\to\infty} \frac{\omega(c)}{\pi(c)} \leq 0$. But now 
\begin{center}
$\lim\limits_{c\to\infty} \Large{\frac{\# \{\varphi_{p,c}(x)\in \mathbb{Z}[x] \ : \ 3\leq p\leq c \ \text{and} \ N_{c}(p) \ = \ 3\}}{\Large{\# \{\varphi_{p,c}(x) \in \mathbb{Z}[x] \ : \ 3\leq p\leq c \}}}} =\lim\limits_{c\to\infty} \frac{\omega(c)}{\pi(c)} \leq 0$.
\end{center}Moreover, we also observe that the number $\# \{\varphi_{p,c}(x)\in \mathbb{Z}[x] : 3\leq p\leq c \ \text{and} \ N_{c}(p) \ = \ 3\}\geq 1$, and so have 
\begin{center}
$\lim\limits_{c\to\infty}\Large{\frac{\# \{\varphi_{p,c}(x) \in \mathbb{Z}[x] \ : \ 3\leq p\leq c \ \text{and} \ N_{c}(p) \ = \ 3\}}{\Large{\# \{\varphi_{p,c}(x) \in \mathbb{Z}[x] \ : \ 3\leq p\leq c \}}}}\geq \lim\limits_{c\to\infty}\frac{1}{\pi(c)} = 0$. But now, we then conclude that the limit 
\end{center}  $\lim\limits_{c\to\infty} \Large{\frac{\# \{\varphi_{p,c}(x) \in \mathbb{Z}[x] \ : \ 3\leq p\leq c \ \text{and} \ N_{c}(p) \ = \ 3\}}{\Large{\# \{\varphi_{p,c}(x) \in \mathbb{Z}[x] \ : \ 3\leq p\leq c \}}}} = 0$ as needed. This then completes the whole proof, as desired.
\end{proof}\noindent Note that we may also certainly interpret Corollary \ref{Thm 4.1} as saying that the probability of choosing randomly a monic polynomial $\varphi_{p,c}(x)$ in the space $\mathbb{Z}[x]$ having exactly three distinct integral fixed points modulo $p$ is zero.
\section{On Densities of Monic Integer Polynomials $\varphi_{p-1,c}(x)$ with $M_{c}(p) = 1$ \& $M_{c}(p) = 2$}\label{sec6}

As in Section \ref{sec5}, motivated by: \say{\textit{What is the density of monics $\varphi_{p-1,c}(x)\in \mathbb{Z}[x]$ having two integral fixed points modulo $p$?}}, we prove in the following corollary that the density of such $\varphi_{p-1,c}(x)\in \mathbb{Z}[x]$ also exists and is zero:
\begin{cor}\label{6.1}
Let $p\geq 5$ be a prime integer. The density of monic polynomials $\varphi_{p-1,c}(x)=x^{p-1}+c\in \mathbb{Z}[x]$ with the number $M_{c}(p) = 2$ exists and is equal to $0\%$ as $c\to \infty$. More precisely, we have 
\begin{center}
    $\lim\limits_{c\to\infty} \Large{\frac{\# \{\varphi_{p-1,c}(x) \in \mathbb{Z}[x]\ : \ 5\leq p\leq c \ and \ M_{c}(p) \ = \ 2\}}{\Large{\# \{\varphi_{p-1,c}(x) \in \mathbb{Z}[x]\ : \ 5\leq p\leq c \}}}} = \ 0.$
\end{center}
\end{cor}
\begin{proof}
By applying a similar argument as in the Proof of Corollary \ref{Thm 4.1}, we then obtain the limit, as desired.
\end{proof} \noindent As before, we note that one may also interpret Corollary \ref{6.1} by saying that the probability of choosing randomly a monic polynomial $\varphi_{p-1,c}(x)$ in the space $\mathbb{Z}[x]$ with exactly two distinct integral fixed points modulo $p$ is zero.

Similarly, the following immediate corollary shows that the probability of choosing randomly a monic integer polynomial $\varphi_{p-1,c}(x)=x^{p-1}+c$ in the space $\mathbb{Z}[x]$ having one integral fixed point modulo $p$ is also zero:

\begin{cor}\label{6.2}
Let $p\geq 5$ be a prime integer. The density of monic polynomials $\varphi_{p-1,c}(x)=x^{p-1}+c\in \mathbb{Z}[x]$ with the number $M_{c}(p) = 1$ exists and is equal to $0\%$ as $c\to \infty$. More precisely, we have 
\begin{center}
    $\lim\limits_{c\to\infty} \Large{\frac{\# \{\varphi_{p-1,c}(x) \in \mathbb{Z}[x]\ : \ 5\leq p\leq c \ and \ M_{c}(p) \ = \ 1\}}{\Large{\# \{\varphi_{p-1,c}(x) \in \mathbb{Z}[x]\ : \ 5\leq p\leq c \}}}} = \ 0.$
\end{center}
\end{cor}
\begin{proof}
As before, $M_{c}(p) = 1$ is as we proved in Theorem \ref{6.0.2}, determined whenever the coefficient $c$ is such that $c-1$ is divisible by any fixed prime $p\geq 5$; and so we may count $\# \{\varphi_{p-1,c}(x) \in \mathbb{Z}[x] : 5\leq p\leq c \ \text{and} \ M_{c}(p) \ = \ 1\}$ by simply counting the number $\# \{\varphi_{p-1,c}(x)\in \mathbb{Z}[x] : 5\leq p\leq c \ \text{and} \ p\mid (c-1) \ \text{for \ any \ fixed} \ c \}$. But now since $c-1<c$, then if the number $\# \{p : 5\leq p\leq c \ \text{and} \ p\mid (c-1) \}< \# \{p : 5\leq p\leq c \ \text{and} \ p\mid c \}$, we then have that
\begin{center}
$\Large{\frac{\# \{\varphi_{p-1,c}(x) \in \mathbb{Z}[x] \ : \ 5\leq p\leq c \ \text{and} \ p\mid (c-1) \ \text{for any fixed} \ c\}}{\Large{\# \{\varphi_{p-1,c}(x) \in \mathbb{Z}[x] \ : \ 5\leq p\leq c \}}}} < \Large{\frac{\# \{\varphi_{p-1,c}(x)\in \mathbb{Z}[x] \ : \ 5\leq p\leq c \ \text{and} \ p\mid c \ \text{for any fixed} \ c \}}{\Large{\# \{\varphi_{p-1,c}(x) \in \mathbb{Z}[x] \ : \ 5\leq p\leq c \}}}}.$ 
\end{center}Letting $c$ tend to infinity on both sides of the above inequality and then applying Corollary \ref{6.1}, we then have 
\begin{center}
$\lim\limits_{c\to\infty}\Large{\frac{\# \{\varphi_{p-1,c}(x) \in \mathbb{Z}[x] \ : \ 5\leq p\leq c \ \text{and} \ p\mid (c-1)\}}{\Large{\# \{\varphi_{p-1,c}(x) \in \mathbb{Z}[x] \ : \ 5\leq p\leq c \}}}} \leq \lim\limits_{c\to\infty}\Large{\frac{\# \{\varphi_{p-1,c}(x)\in \mathbb{Z}[x] \ : \ 5\leq p\leq c \ \text{and} \ p\mid c  \}}{\Large{\# \{\varphi_{p-1,c}(x) \in \mathbb{Z}[x] \ : \ 5\leq p\leq c \}}}}$ = 0; 
\end{center} from which it then follows that 
\begin{center}
    $\lim\limits_{c\to\infty} \Large{\frac{\# \{\varphi_{p-1,c}(x) \in \mathbb{Z}[x]\ : \ 5\leq p\leq c \ \text{and} \ M_{c}(p) \ = \ 1\}}{\Large{\# \{\varphi_{p-1,c}(x) \in \mathbb{Z}[x]\ : \ 5\leq p\leq c \}}}} = \ 0$, and hence showing the limit as desired in this case.
\end{center}Otherwise, if the number $\# \{p : 5\leq p\leq c \ \text{and} \ p\mid c \}< \# \{p : 5\leq p\leq c \ \text{and} \ p\mid (c-1) \}$, we then have that
\begin{center}
$\Large{\frac{\# \{\varphi_{p-1,c}(x) \in \mathbb{Z}[x] \ : \ 5\leq p\leq c \ \text{and} \ p\mid c \ \text{for any fixed} \ c\}}{\Large{\# \{\varphi_{p-1,c}(x) \in \mathbb{Z}[x] \ : \ 5\leq p\leq c \}}}} < \Large{\frac{\# \{\varphi_{p-1,c}(x)\in \mathbb{Z}[x] \ : \ 5\leq p\leq c \ \text{and} \ p\mid (c-1) \ \text{for any fixed} \ c \}}{\Large{\# \{\varphi_{p-1,c}(x) \in \mathbb{Z}[x] \ : \ 5\leq p\leq c \}}}}.$
\end{center}So now, taking limit as $c\to \infty$ on both sides of the above inequality and applying Corollary \ref{6.1} and then applying a similar argument as in the Proof of Corollary \ref{Thm 4.1} to obtain an upper bound zero, we then obtain 
\begin{center}
$\lim\limits_{c\to\infty}\Large{\frac{\# \{\varphi_{p-1,c}(x) \in \mathbb{Z}[x] \ : \ 5\leq p\leq c \ \text{and} \ p\mid (c-1)\}}{\Large{\# \{\varphi_{p-1,c}(x) \in \mathbb{Z}[x] \ : \ 5\leq p\leq c \}}}} = 0$ as needed. This then completes the whole proof, as desired.
\end{center}  
\end{proof}

\section{On Density of Integer Monics $\varphi_{p,c}(x)$ with $N_{c}(p) = 0$ and $\varphi_{p-1,c}(x)$ with $M_{c}(p) = 0$}\label{sec7}

Recall in Corollary \ref{Thm 4.1} that a density of $0\%$ of monic integer polynomials $\varphi_{p,c}(x)$ have exactly three integral fixed points modulo $p$; and so the density of monic integer polynomials $\varphi_{p,c}(x)-x$ that are reducible modulo $p$ is $0\%$. So now, we also wish to determine: \say{\textit{What is the density of monic integer polynomials $\varphi_{p,c}(x)=x^p + c$ with no integral fixed points modulo $p$?}} The following corollary shows that the probability of choosing randomly a monic integer polynomial $\varphi_{p,c}(x)=x^p + c$ such that $\mathbb{Q}[x]\slash (\varphi_{p, c}(x)-x)$ is a number field of odd degree $p$ is $1$: 
\begin{cor}\label{8.1}
Let $p\geq 3$ be a prime integer. Then the density of monic polynomials $\varphi_{p,c}(x) = x^p + c\in \mathbb{Z}[x]$ with the number $N_{c}(p) = 0$ exists and is equal to $100 \%$ as $c\to \infty$. More precisely, we have 
\begin{center}
    $\lim\limits_{c\to\infty} \Large{\frac{\# \{\varphi_{p,c}(x)\in \mathbb{Z}[x] \ : \ 3\leq p\leq c \ and \ N_{c}(p) \ = \ 0 \}}{\Large{\# \{\varphi_{p,c}(x) \in \mathbb{Z}[x] \ : \ 3\leq p\leq c \}}}} = \ 1.$
\end{center}
\end{cor}
\begin{proof}
Since the number $N_{c}(p) = 3$ or $0$ for any given prime integer $p\geq 3$ and since we also proved the density in Corollary \ref{Thm 4.1}, we then obtain the desired density (i.e., we obtain that the limit exists and is equal to 1). 
\end{proof}
\noindent Note that the foregoing corollary also shows that there are infinitely many polynomials $\varphi_{p,c}(x)$ over $\mathbb{Q}$ such that for $f(x) = \varphi_{p,c}(x)-x = x^p-x+c$, the induced quotient ring $K_{f} = \mathbb{Q}[x]\slash (f(x))$ is an algebraic number field of odd degree $p$. Comparing the densities in Corollaries \ref{Thm 4.1} and \ref{8.1}, we may also observe that in the whole family of monic integer polynomials $\varphi_{p,c}(x) = x^p +c$, almost all such monics have no integral fixed points modulo $p$ (i.e., have no rational roots); and so almost all such monic polynomials $f(x)$ are irreducible over $\mathbb{Q}$. Consequently, this may also imply that the average value of $N_{c}(p)$ in the whole family of polynomials $\varphi_{p,c}(x)\in \mathbb{Z}[x]$ is zero.

Similarly, recall in Corollary \ref{6.1} and \ref{6.2} that a density of $0\%$ of monic integer polynomials $\varphi_{p-1,c}(x)$ have $M_{c}(p) = 2$ or $1$, respectively; and so the density of polynomials $\varphi_{p-1, c}(x)-x\in \mathbb{Z}[x]$ that are reducible modulo $p$ is $0\%$. So now, we also wish to determine: \say{\textit{What is the density of monic polynomials $\varphi_{p-1,c}(x)\in \mathbb{Z}[x]$ with no integral fixed points modulo $p$?}} The following corollary shows that the probability of choosing randomly a monic polynomial $\varphi_{p-1,c}(x)\in \mathbb{Z}[x]$ so that $\mathbb{Q}[x]\slash (\varphi_{p-1, c}(x)-x)$ is a number field of even degree $p-1$ is also $1$:
\begin{cor} \label{5.1}
Let $p\geq 5$ be a prime integer. The density of monic polynomials $\varphi_{p-1, c}(x) = x^{p-1}+c\in \mathbb{Z}[x]$ with the number $M_{c}(p) = 0$ exists and is equal to $100 \%$ as $c\to \infty$. More precisely, we have 
\begin{center}
    $\lim\limits_{c\to\infty} \Large{\frac{\# \{\varphi_{p-1, c}(x)\in \mathbb{Z}[x] \ : \ 5\leq p\leq c \ and \ M_{c}(p) \ = \ 0 \}}{\Large{\# \{\varphi_{p-1,c}(x) \in \mathbb{Z}[x] \ : \ 5\leq p\leq c \}}}} = \ 1.$
\end{center}
\end{cor}
\begin{proof}
Recall that the number $M_{c}(p) = 1, 2$ or $0$ for any given prime $p\geq 5$ and since we also proved the densities in Cor. \ref{6.1} and \ref{6.2}, we now obtain the desired density (i.e., we get that the limit exists and is equal to 1).
\end{proof}
\noindent As before, Corollary \ref{5.1} also shows that there are infinitely many monic polynomials $\varphi_{p-1,c}(x)$ over $\mathbb{Q}$ such that for $g(x) = \varphi_{p-1,c}(x)-x = x^{p-1}-x+c$, the induced quotient ring $L_{g} = \mathbb{Q}[x]\slash (g(x))$ is an algebraic number field of even degree $p-1$. Again, if we compare the densities in Cor. \ref{6.1}, \ref{6.2} and \ref{5.1}, we may then see that in the whole family of monic  polynomials $\varphi_{p-1,c}(x) = x^{p-1} +c\in \mathbb{Z}[x]$, almost all such monics have no integral fixed points modulo $p$ (i.e., have no $\mathbb{Q}$-roots); and so almost all monic polynomials $g(x)$ are irreducible over $\mathbb{Q}$. But this may also imply that the average value of $M_{c}(p)$ in the whole family of monics $\varphi_{p-1,c}(x)$ is also zero.

As always a central theme in algebraic number theory that whenever one is studying an algebraic number field $K$ of some interest, one must simultaneously try to describe very precisely what the associated ring $\mathcal{O}_{K}$ of integers is; and this is because this ring of integers $\mathcal{O}_{K}$ is classically known to describe naturally the arithmetic of the underlying number field $K$. However, accessing $\mathcal{O}_{K}$ in practice from a computational point of view is a well-known extremely involved problem. So now, in our case here, it then follows $K_{f}$ has a ring of integers $\mathcal{O}_{K_{f}}$, and moreover by Bhargava-Shankar-Wang \cite{sch1}, we then also obtain the following corollary showing the probability of choosing randomly a polynomial $f\in \mathbb{Z}[x]$ arising from a polynomial discrete dynamical system in Section \ref{sec2} (and ascertained by Corollary \ref{8.1}), such that $\mathbb{Z}[x]\slash (f(x))$ is the ring of integers of $K_{f}$, is $\approx 60.7927\%$:
\begin{cor} \label{8.2}
Assume Corollary \ref{8.1}. When integer polynomials $f(x)$ are ordered by height $H(f) = |c|^{1\slash p}$ as defined in \textnormal{\cite{sch1}}, the density of such polynomials $f(x)$ such that $\mathbb{Z}[x]\slash (f(x))$ is the ring of integers of $K_{f}$ is $\zeta(2)^{-1}$. 
\end{cor}

\begin{proof}
Since from Corollary \ref{8.1} we know that there are infinitely many monic polynomials $f(x)$ over $\mathbb{Z}$ (and hence over $\mathbb{Q}$) such that $K_{f} = \mathbb{Q}[x]\slash (f(x))$ is an algebraic number field of degree deg$(f)$; and moreover associated to $K_{f}$ is the ring $\mathcal{O}_{K_{f}}$ of integers. Now applying a remarkable result of Bhargava-Shankar-Wang [\cite{sch1}, Theorem 1.2] to the underlying family of monic integer polynomials $f(x)$ ordered by height $H(f) = |c|^{1\slash p}$ such that $\mathcal{O}_{K_{f}} = \mathbb{Z}[x]\slash (f(x))$, we then obtain that the density of such polynomials $f$ is $\zeta(2)^{-1} \approx 60.7927\%$, as needed. 
\end{proof}

As with $K_{f}$, every number field $L_{g}$ induced by a polynomial $g$, is naturally equipped with the ring of integers $\mathcal{O}_{L_{g}}$; and moreover this ring $\mathcal{O}_{L_{g}}$ may also be difficult to compute in practice. So now, as before we take great advantage of a density result in [\cite{sch1}, Theorem 1.2] and then obtain the following corollary showing the probability of choosing randomly $g\in \mathbb{Z}[x]$ arising from a polynomial discrete dynamical system in Section \ref{sec3} (and also ascertained by Corollary \ref{5.1}), such that $\mathbb{Z}[x]\slash (g(x))$ is the ring of integers of $L_{g}$, is also $\approx 60.7927\%$:

\begin{cor}
Assume Corollary \ref{5.1}. When integer polynomials $g(x)$ are ordered by height $H(g) = |c|^{1\slash (p-1)}$ as defined in \textnormal{\cite{sch1}}, the density of polynomials $g(x)$ such that $\mathbb{Z}[x]\slash (g(x))$ is the ring of integers of $L_{g}$ is $\zeta(2)^{-1}$. 
\end{cor}

\begin{proof}
By applying a similar argument as in the Proof of Corollary \ref{8.2}, we then obtain density, as needed.
\end{proof}

\section{On the Number of Number fields $K_{f}$ \& $L_{g}$ with Bounded Absolute Discriminant}\label{sec8}

Recall from Corollary \ref{8.1} that there is an infinite family of irreducible monic integer polynomials $f(x) = x^p-x + c$ such that the quotient ring $K_{f} = \mathbb{Q}[x]\slash (f(x))$ induced by $f$ is an algebraic number field of odd prime degree $p$. Moreover, we may also recall from Corollary \ref{5.1} that one can always find an infinite family of irreducible monic integer polynomials $g(x) = x^{p-1}-x + c$ such that the field extension $L_{g}$ over $\mathbb{Q}$ arising from $g$ is an algebraic number field of even degree $p-1\geq 4$. In this section, we wish to study the problem of counting number fields; a problem that is originally from and is of very serious interest in arithmetic statistics. Inspired by Bhargava-Shankar-Wang (BSW) \cite{sch}, we then wish to count here the number of primitive number fields $K_{f}$ induced by irreducible polynomials $f\in \mathbb{Z}[x]$ arising from a polynomial discrete dynamical system in Section \ref{sec2} (and ascertained by Corollary \ref{8.1}), with bounded absolute discriminant. With that in mind, we then obtain:

\begin{cor}\label{cor6}
Assume Corollary \ref{8.1}, and let $K_{f} = \mathbb{Q}[x]\slash (f(x))$ be a primitive number field with discriminant $\Delta(K_{f})$. Then up to isomorphism classes of number fields, we have that $\# \{ K_{f} \ : \ |\Delta(K_{f})| < X \} \ll  X^{p\slash(2p-2)}$. 
\end{cor}

\begin{proof}
From Corollary \ref{8.1}, we know that there are infinitely many monics $f(x) = x^p - x + c$ over $\mathbb{Z}$ (and so over $\mathbb{Q}$) such that $K_{f}$ is a number field of odd prime degree $p\geq 3$. But now when the underlying number fields $K_{f}$ are primitive, we may then apply an argument of Bhargava-Shankar-Wang [\cite{sch}, Page 2] to show that up to isomorphism classes of number fields, the total number of such primitive number fields $K_{f}$ with $|\Delta(K_{f})| < X$ is bounded above by the number of monic integer polynomials $f$ of degree $p\geq 3$ with height $H(f) \ll X^{1\slash(2p-2)}$ and vanishing subleading coefficient. Now since $H(f) = |c|^{1\slash p}$ and so $|c| \ll X^{p\slash(2p-2)}$, we then have that the size $\# \{f(x)\in \mathbb{Z}[x] \ : H(f) \ll X^{1\slash(2p-2)} \} = \# \{f(x)\in \mathbb{Z}[x] \ : |c| \ll X^{p\slash(2p-2)} \} \ll X^{p\slash(2p-2)}$. Hence, we obtain that the number $\# \{ K_{f} : |\Delta(K_{f})| < X \} \ll  X^{p\slash(2p-2)}$ up to isomorphism classes of number fields, as needed.
\end{proof}

By applying a similar argument as in Corollary \ref{cor6}, we then also obtain the following corollary on the number of primitive number fields $L_{g}$ induced by irreducible polynomials $g(x)\in \mathbb{Z}[x]$ arising from a polynomial discrete dynamical system in Section \ref{sec3} (and ascertained by Corollary \ref{5.1}), with bounded absolute discriminant:

\begin{cor}
Assume Corollary \ref{5.1}, and let $L_{g} = \mathbb{Q}[x]\slash (g(x))$ be a primitive number field with discriminant $\Delta(L_{g})$. Then up to isomorphism classes of number fields, we have that $\# \{ L_{g} \ : \ |\Delta(L_{g})| < X \} \ll  X^{(p-1)\slash(2p-4)}$. 
\end{cor}

We recall in algebraic number theory that an algebraic number field $K$ is called \say{\textit{monogenic}} if there exists an algebraic number $\alpha \in K$ such that the ring $\mathcal{O}_{K}$ of integers is the subring $\mathbb{Z}[\alpha]$ generated by $\alpha$ over $\mathbb{Z}$, i.e., $\mathcal{O}_{K}= \mathbb{Z}[\alpha]$. In Corollary \ref{cor6}, we counted primitive number fields $K_{f}$ (i.e., $K_{f} = \mathbb{Q}(\alpha)$ for some $\alpha \in K$) with $|\Delta(K_{f})| < X$. So now, we wish to count the number of number fields $K_{f}$ induced by irreducible polynomials $f\in \mathbb{Z}[x]$ arising from a polynomial discrete dynamical system in Section \ref{sec2} (and ascertained by Corollary \ref{8.1}), that are monogenic with $|\Delta(K_{f})| < X$ and with Galois group Gal$(K_{f}\slash \mathbb{Q})$ equal to the symmetric group $S_{p}$. A hard counting task that we very easily tackle by taking great advantage of a result of (BSW)[\cite{sch1}, Corollary 1.3]:

\begin{cor}\label{8.3}
Assume Corollary \ref{8.1}. Then the number of isomorphism classes of algebraic number fields $K_{f}$ of degree $p\geq 3$ and with $|\Delta(K_{f})| < X$ that are monogenic and have associated Galois group $S_{p}$ is $\gg X^{\frac{1}{2} + \frac{1}{p}}$.
\end{cor}

\begin{proof}
By Corollary \ref{8.1}, it then also follows that there are infinitely many polynomials $f\in \mathbb{Z}[x]$ such that $K_{f}$ is a degree-$p$ number field. This also means that the set of fields $K_{f}$ of degree $p$ is not empty. But now applying [\cite{sch1}, Corollary 1.3] to the underlying set of fields $K_{f}$ with $|\Delta(K_{f})| < X$ that are monogenic and with Galois group $S_{p}$, we then obtain that the number of isomorphism classes of such fields $K_{f}$ is $\gg X^{\frac{1}{2} + \frac{1}{p}}$, as needed.
\end{proof}

Similarly, by again we again taking great advantage of that same result of (BSW)[\cite{sch1}, Corollary 1.3], we then also obtain the following corollary on the number of number fields $L_{g}$ induced by irreducible polynomials $g\in \mathbb{Z}[x]$ arising from a polynomial discrete dynamical system in Section \ref{sec3} (and ascertained by Corollary \ref{5.1}), that are monogenic with $|\Delta(L_{g})| < X$ and such that Galois group Gal$(L_{g}\slash \mathbb{Q})$ is equal to symmetric group $S_{p-1}$: 

\begin{cor}
Assume Corollary \ref{5.1}. Then the number of isomorphism classes of algebraic number fields $L_{g}$ of degree $p-1\geq 4$ and $|\Delta(L_{g})| < X$ that are monogenic and have associated Galois group $S_{p-1}$ is $\gg X^{\frac{1}{2} + \frac{1}{p-1}}$.
\end{cor}

\begin{proof}
By applying a similar argument as in the Proof of Corollary \ref{8.3}, we then obtain the count, as needed.
\end{proof}

\addcontentsline{toc}{section}{Acknowledgements}
\section*{\textbf{Acknowledgements}}
I'm deeply indebted to my long-time great advisors, Dr. Ilia Binder and Dr. Arul Shankar, for all their boundless generosity, friendship and for all the inspiring weekly conversations and along with Dr. Jacob Tsimerman for always very uncomprimisingly supporting my professional and philosophical-mathematical research endeavours. I'm very grateful to Dr. Shankar for bringing my attention to a Number theory course (MAT1200HS, 2023/24) and for strongly suggesting it to me. This suggestion would've never easily be realized if the course instructor at that time had not been unsurpassably generous in terms of allowing anyone to attend his class and in terms of his friendliness; and so I'm once again deeply indebted to Dr. Shankar for everything and to that instructor, namely, Dr. Tsimerman for his very enlightening Algebraic Number Theory class and conversations. Not to forget, I am very grateful to the very vibrant Dept. of Mathematical and Computational Sciences (MCS) at the University of Toronto, Mississuaga; and in particular, I'm very grateful and deeply indebted to Dr. Yael Karshon (for the truly invaluable advice about the state of professional mathematics), Dr. Ilia Binder, Dr. Arul Shankar, Dr. Marina Tvalavadze, Dr. Alex Rennet, Dr. Michael Gr\"{o}echenig, Dr. Julie Desjardins, Dr. Duncan Dauvergne, Dr. Ke Zhang, Dr. Jaimal Thind, for everything. I'm truly very grateful to the Office of the President at the UofT for the amazing hospitality during my in-person visit to converse with President Meric S. Gertler on \say{Leadership and Duty}. I'm very grateful to Dr. Meric for the very enlightening great conversations on the insurmountable importance of collaboration, inclusion and diversity in educational settings, and also more importantly for not only envisioning with Rose M. Patten for more representation of minority groups in U of T Science and Mathematics departments but also for taking very serious practical steps toward realizing such a vision with integrity, as I've thoroughly witnessed during the fantastic and progressive chairship of Dr. Binder with his team in Dept. of MCS. Lastly, I’m very grateful and deeply indebted to the former U of T Mississauga Registrar and Director of Enrolment Management, namely, Lorretta Neebar, for everything. Last but not the least, I'm truly very grateful and indebted to Dr. Michael Bumby for the great life-philosophical conversations, and more importantly for being a truly an amazing life coach and a great friend.
As a graduate research student, this work and my studies are hugely and wholeheartedly funded by Dr. Binder and Dr. Shankar. As part of the first harvest of an upcoming long harvest, I very happily dedicate this article to the Department of Mathematical and Computational Sciences (MCS) and the Presidency of Dr. Meric S. Gerlter at the University of Toronto! Any opinions expressed in this article belong solely to me, the author, Brian Kintu; and should never be taken at all as a reflection of the views of anyone that’s been happily acknowledged by me.

\bibliography{References}
\bibliographystyle{plain}

\noindent Dept. of Math. and Comp. Sciences (MCS), University of Toronto, Mississauga, Canada \newline
\textit{E-mail address:} \textbf{brian.kintu@mail.utoronto.ca}\newline 
\date{\small{\textit{January 15, 2026.}}}

\end{document}